\documentclass[twoside]{article}

\usepackage{RR}
\usepackage{hyperref}
\usepackage{graphicx}
\usepackage{epsfig}
\usepackage{subfigure}
\usepackage{amsmath,amsgen,amstext,amssymb}
\usepackage{amsfonts}
\usepackage{latexsym}

\RRdate{June 2012}

\RRauthor{
Konstantin Avrachenkov\thanks{INRIA Sophia Antipolis, France, K.Avrachenkov@sophia.inria.fr}%
\and
Alexey Piunovskiy\thanks{The University of Liverpool, Department of Maths Sciences,
M\&O Building, L69 7ZL, UK, piunov@liverpool.ac.uk}
\and
Yi Zhang\thanks{The University of Liverpool, Department of Maths Sciences,
M\&O Building, L69 7ZL, UK, yi.zhang@liverpool.ac.uk}
}

\authorhead{K. Avrachenkov \& A. Piunovskiy \& Y. Zhang}

\RRtitle{Les Processus de Markov avec Red\'emarrage}

\RRetitle{Markov Processes with Restart}

\titlehead{Markov Processes with Restart}

\RRresume{Nous consid\'erons un processus de Markov en temps continu qui est g\'en\'eral, honn\^ete et homog\`ene.
Le  processus est forc\'e \`a red\'emarrer \`a partir d'une distribution donn\'ee \`a des moments de temps
g\'en\'er\'es par un processus ind\'ependant de Poisson. La motivation pour \'etudier ce type de processus
vient de la mod\'elisation de la mobilit\'e humaine et animale, des proc\'edures avec le red\'emarrage dans
les protocoles de communication et de l'application des marches al\'eatoires avec le red\'emarrage
en recherche d'information. Nous fournissons une connexion entre les fonctions de probabilit\'e de transition
du processus de Markov originale et du processus modifi\'e avec le red\'emarrage. Nous donnons des
expressions explicites de la mesure invariante du processus modifi\'e. Lorsque le processus \'evolue dans
l'espace Euclidien, il y a \'egalement une expression explicite pour les moments du processus modifi\'e.
Nous montrons que le processus modifi\'e est \`a la fois Harris positif r\'ecurrent et ergodique exponentiel
avec l'indice \'egal \`a (ou plus grand que) le taux de red\'emarrage. Enfin, nous illustrons les r\'esultats g\'en\'eraux
par les cas classiques de mouvements Browniens standard et g\'eom\'etrique.}

\RRabstract{We consider a general honest homogeneous continuous-time Markov process with restarts.
The process is forced to restart from a given distribution at time moments generated by
an independent Poisson process. The motivation to study such processes comes from modeling
human and animal mobility patterns, restart processes in communication protocols, and
from application of restarting random walks in information retrieval. We provide a connection
between the transition probability functions of the original Markov process and the modified
process with restarts. We give closed-form expressions for the invariant probability
measure of the modified process. When the process evolves on the Euclidean space there is
also a closed-form expression for the moments of the modified process. We show that the
modified process is always positive Harris recurrent and exponentially ergodic with the index
equal to (or bigger than) the rate of restarts. Finally, we illustrate the general results
by the standard and geometric Brownian motions.}

\RRmotcle{Processus de Markov avec Red\'emarrage, Positive Harris Recurrence, Ergodicit\'e
Exponentiel, Mouvements Browniens Standard et G\'eom\'etrique}

\RRkeyword{Markov Processes with Restart, Positive Harris Recurrence, Exponential Ergodicity,
Standard and Geometric Brownian Motions}

\RRprojet{Maestro}

\RCSophia

\begin{document}

\RRNo{8000}
\makeRR 

\newcommand{\beq}{\begin{equation}}
\newcommand{\deq}{\end{equation}}

\newcommand{\baq}{\begin{eqnarray}}
\newcommand{\daq}{\end{eqnarray}}

\newcommand{\baqm}{\begin{eqnarray*}}
\newcommand{\daqm}{\end{eqnarray*}}

\newcommand{\up}[1]{\overline{#1}}
\newcommand{\down}[1]{\underline{#1}}

\newcommand{\ti}[1]{\tilde{#1}}
\newcommand {\R} {\rm I\!R}
\newcommand{\qed}{\rightline{$\Box$}}
\newcommand{\dfn}{\stackrel{\Delta}{=}}

\newcommand{\eps}{\varepsilon}
\newenvironment{pf}{{\bf Proof. }}{\hfill $\square$\medskip}
\newtheorem{thm}{Theorem}
\newtheorem{ass}{Assumption}
\newtheorem{defi}{Definition}
\newtheorem{prop}{Proposition}
\newtheorem{cor}{Corollary}
\newtheorem{lem}{Lemma}
\newtheorem{rem}{Remark}

\section{Introduction}

Many phenomena in nature and technology can be modeled by Markov processes which restart
from time to time. Human and animal movements can be modeled by Markov processes that
restart from some locations \cite{GHB08,WBC10}. Facing the congestion of the Internet traffic,
the Internet users tend to restart their sessions \cite{KR01,MH01} or a protocol governing the
rate of information transmission invokes restart routines (i.e., Slow-Start routine in the Internet
Transmission Control Protocol) \cite{S94}. The celebrated PageRank algorithm \cite{BP98} in
information retrieval models the behavior of a random surfer, who decides to restart web surfing
from time to time.
The heat kernel PageRank \cite{C07} is a continuous-time analog of the original discrete-time algorithm.
The restart policy is also used to speedup the Las Vegas type randomized
algorithms \cite{Alt96,LSZ93}.

Motivated by the above mentioned phenomena, we study a general honest homogeneous continuous-time
Markov process with restarts. We assume that the process of restarts is modeled by
an independent Poisson process. In the next Section~2 we derive a formula which makes a connection
between the transition probability functions of the original Markov process and the Markov process
modified by restarts. As a corollary we obtain a closed-form expression for the invariant probability
measure of the modified process. We also show that the modified process is always positive Harris
recurrent and exponentially ergodic with the index equal to (or bigger than) the rate of restarts.
In Section~3 we provide
bounds for the moments of the modified process and investigate the limiting behaviour of the
modified process when the rate of restarts goes to zero. We conclude the paper with Section~4 where
we consider the restart modifications of the two classical Markov processes: the standard Brownian motion
with drift and the geometric Brownian motion. It is very curious to observe that in the case of the geometric
Brownian motion the number of finite moments of the modified process depends on the rate of restarts.
Thus, even a small change in the value of the restart rate can lead to dramatic changes of the properties
of the geometric Brownian motion with restarts.

\section{Main results}

In accordance with \cite{Kuz80}, we consider a homogeneous continuous-time Markov process $\{X(t) : t \in [0,\infty)\}$ in a Borel space $(E,{\cal E})$ characterized by the initial distribution $\gamma(\cdot)$, and the honest transition (probability) function $P(t,x,dy),$
satisfying the following properties:
\begin{enumerate}
\item $P(t,x,\cdot)$ is a probability measure on $E$,
\item $P(0,x,\Gamma)=1\{x \in \Gamma\}$,
\item for each fixed $\Gamma\in {\cal E}$, $t\in[0,\infty),$ $P(t,x,\Gamma)$ is a jointly measurable function with
respect to $(t,x)\in[0,\infty)\times E$,
\item and the Chapman-Kolmogorov equation takes place
$$
P(t+s,x,\Gamma) = \int_E P(s,y,\Gamma) P(t,x,dy).
$$
\end{enumerate}

Throughout this paper, all the processes are from the same probability space $(\Omega,{\cal F},P)$. We further use the notations $P_x,E_x,Var_x$ when the initial state of the concerned process is $x\in E.$

The goal of the present work is to analyse a modification of the Markov process $\{X(t) : t \in [0,\infty)\}$ introduced in the above by forcing
the process to restart with a given restart distribution $\nu(\Gamma)$ after an exponentially distributed
random time, i.e., the process counting restarts represents a standard Poisson process with the rate, say, $\lambda>0$, independent of $\{X(t) : t \in [0,\infty)\}$.  In the following theorem we characterize the transition function of the modified Markov process, which is denoted as $\{\tilde{X}(t) : t \in [0,\infty)\}.$


\begin{thm}\label{thm:Ptilde}
Let $P(t,x,\Gamma)$ be the transition function of the Markov process $\{X(t) : t \in [0,\infty)\}$ taking values in a Borel space
$(E,{\cal E})$. Then, the modified Markov process $\{\tilde{X}(t),t\in [0,\infty)\}$ that restarts from a distribution $\nu$ after (independent) exponentially distributed random times with mean $1/\lambda$ has the following transition function
\begin{equation}\label{eq:Ptilde}
\tilde P_\nu(t,x,\Gamma) = e^{-\lambda t} P(t,x,\Gamma)
+\int_E \int_0^t \lambda e^{-\lambda s} P(s,y,\Gamma) \, ds \, \nu(dy).
\end{equation}
\end{thm}
{\bf Proof:} Denote by $N(t)$ the number of restarts up to time $t$ and by $S$ the time elapsed since the time moment
of the last restart before $t$. Then, we have
$$
\tilde P_\nu(t,x,\Gamma) = P[\{\tilde X(t) \in \Gamma\} \cap \{N(t)=0\}|\tilde X(0)=x]
$$
\begin{equation}\label{eq:Ptildedecomp}
+ P[\{\tilde X(t) \in \Gamma\} \cap \{N(t)>0\}|\tilde X(0)=x].
\end{equation}
The first term in the above equation is given by
$$
P[\{\tilde X(t) \in \Gamma\} \cap \{N(t)=0\}|\tilde X(0)=x] = e^{-\lambda t}P(t,x,\Gamma).
$$
Let us now calculate the second term in (\ref{eq:Ptildedecomp}). Note that the second term in (\ref{eq:Ptildedecomp})
can be expressed as
\begin{equation}\label{eq:Ptilde2term}
P[\{\tilde X(t) \in \Gamma\} \cap \{N(t)>0\}|\tilde X(0)=x]=
\int_E \int_0^t P(s,y,\Gamma) \, dF(s) \, \nu(dy),
\end{equation}
where $F(s)=P[S \le s|N(t)>0]$. Let $T_k$ denote the $k$-th restart moment from time zero.
If $N(t)=n>0$, then the restart times $T_1,...,T_n$ have the same distribution as the order
statistics corresponding to $n$ independent random variables uniformly distributed on the
interval $(0,t)$ \cite{R96}, and therefore,
$$
P[T_n \le \tau | N(t)=n] = \left(\frac{\tau}{t}\right)^n,~\tau\in [0,t],~n=1,2,\dots.
$$
Then,
$$
P[S \le s | N(t)=n] = P[T_n > t-s | N(t)=n] = 1-\left(\frac{t-s}{t}\right)^n, ~n=1,2,\dots,
$$
and, consequently,
$$
F(s) = \sum_{n=1}^\infty \frac{(\lambda t)^n}{n!} e^{-\lambda t} \left[1-\left(\frac{t-s}{t}\right)^n\right]
=e^{-\lambda t}[e^{\lambda t}-e^{\lambda (t-s)}]=1-e^{-\lambda s}.
$$
Hence,
$$
P[\{\tilde X(t) \in \Gamma\} \cap \{N(t)>0\}|\tilde X(0)=x] =
\int_E \int_0^t \lambda e^{-\lambda s} P(s,y,\Gamma) \, ds \, \nu(dy),
$$
and the formula (\ref{eq:Ptilde}) follows.

\qed


\begin{cor}\label{cor:q} The measure
\begin{equation}\label{eq:qdef}
q_\nu(\Gamma) = \int_E \int_0^\infty \lambda e^{-\lambda s} P(s,y,\Gamma) \, ds \, \nu(dy)
\end{equation}
is an invariant probability measure for $\tilde P_\nu(t,x,\Gamma)$.
\end{cor}
{\bf Proof:} Since $P(s,y,\cdot)$ is a probability measure for all $s$ and $y$, by the Vitali-Hahn-Saks theorem $q_\nu(\cdot)$
is a probability measure as well. Let us show that it is indeed an invariant measure.
It is enough to show this for $\nu(dy)=1\{w \in dy\}$, where $w\in E$ is fixed, and in that case we also denote  $q_\nu$ by $q_w$ and $\tilde{P}_{\nu}$ by $\tilde{P}_w$ for brevity. By Theorem \ref{thm:Ptilde}, we obtain
$$
\int_E q_w(dz) \tilde P_{w}(t,z,\Gamma)=
\int_E q_w(dz) \left[e^{-\lambda t} P(t,z,\Gamma)+\int_0^t \lambda e^{-\lambda s} P(s,w,\Gamma) \, ds\right]
$$
$$
=e^{-\lambda t} \int_E q_w(dz) P(t,z,\Gamma) + \int_0^t \lambda e^{-\lambda s} P(s,w,\Gamma) \, ds.
$$
Substituting into the above equation
$$
q_w(\Gamma) = \int_0^\infty \lambda e^{-\lambda s} P(s,w,\Gamma) \, ds
$$
and using the Chapman-Kolmogorov equation  yields
$$
\int_E q_w(dz) \tilde P_w(t,z,\Gamma)=
e^{-\lambda t} \int_0^\infty \lambda e^{-\lambda s} P(s+t,w,\Gamma) \, ds
+ \int_0^t \lambda e^{-\lambda s} P(s,w,\Gamma) \, ds
$$
$$
=\int_t^\infty \lambda e^{-\lambda s'} P(s',w,\Gamma) \, ds'
+ \int_0^t \lambda e^{-\lambda s} P(s,w,\Gamma) \, ds
$$
$$
=\int_0^\infty \lambda e^{-\lambda s} P(s,w,\Gamma) \, ds = q_w(\Gamma),
$$
which concludes the proof of the corollary.

\qed

\begin{rem}
We note that alternatively formula (\ref{eq:qdef}) can be rewritten as
$$
q_\nu(\Gamma) = \int_E \lambda R(y,\Gamma) \, \nu(dy),
$$
where $R(y,\Gamma)=\int_0^\infty e^{-\lambda s} P(s,y,\Gamma) \, ds$ is the resolvent operator.
\end{rem}

Now we are ready to prove that the modified process is positive Harris recurrent
and exponentially ergodic with index $\lambda$. Before this, let us remind the
definitions.

Recall that a homogeneous continuous-time Markov process $\{\tilde{X}(t),t\in [0,\infty)\}$ is called Harris recurrent
if there exists a non-trivial $\sigma$-finite (recurrence) measure $\mu$ on ${\cal E}$ such that
for each $x\in E,$
$$
\tau_\Gamma:=\inf\{t\ge 0: \tilde{X}(t)\in \Gamma\}<\infty, \quad P_x-\mbox{a.s.},
$$
whenever $\mu(\Gamma)>0$, see \cite{Glynn:1994,MT93}.

Then, if a Harris recurrent process admits an invariant probability measure, it must be unique,
and the process is further called positive Harris recurrent, see \cite{Glynn:1994,MT93}.

A homogeneous continuous-time Markov process $\{\tilde{X}(t),t\in [0,\infty)\}$ is called exponentially ergodic
with index $\alpha$ if there exist a probability measure $\mu(\cdot)$ on ${\cal E}$,
a finite valued function $M(\cdot)$ on E, and a constant $\alpha > 0$, satisfying
$||\tilde{P}(t,x,\cdot)-\mu(\cdot)||_{TV}\le M(x)e^{-\alpha t}$ for every $x\in E,$ where $||\cdot||_{TV}$
denotes the total variation norm, see \cite{Down:1995}.

\begin{thm}
The modified Markov process $\{\tilde{X}(t),t\in [0,\infty)\}$ is positive Harris recurrent
and exponentially ergodic with index equal to $($or bigger than$)$ the rate of restarts $\lambda,$
and the following inequality takes place
\begin{equation}\label{eq:uniformineq}
|q_\nu(\Gamma)-\tilde P_\nu(t,x,\Gamma)| \le e^{-\lambda t},
\quad \forall x \in E, \ \forall \Gamma \in {\cal E}.
\end{equation}
\end{thm}
{\bf Proof:}
The modified process $\{\tilde{X}(t),t\in [0,\infty)\}$ under consideration is Harris recurrent with the recurrence measure given by the restart distribution $\mu(\cdot)=\nu(\cdot)$. From Corollary~1 we conclude that
$\{\tilde{X}(t),t\in [0,\infty)\}$ is positive Harris recurrent with the unique invariant probability measure
given by (\ref{eq:qdef}). Let us now prove (\ref{eq:uniformineq}). We see
$$
|q_\nu(\Gamma)-\tilde P_\nu(t,x,\Gamma)|=
|\int_E \int_t^\infty \lambda e^{-\lambda s} P(s,y,\Gamma) \, ds \, \nu(dy)-e^{-\lambda t} P(t,x,\Gamma)|
$$
$$
=|\int_E \int_t^\infty \lambda e^{-\lambda s} P(s,y,\Gamma) \, ds \, \nu(dy)
-\int_t^\infty \lambda e^{-\lambda s}P(t,x,\Gamma) \, ds|
$$
$$
=|\int_t^\infty \lambda e^{-\lambda s} \left[ \int_E P(s,y,\Gamma) \, \nu(dy)-P(t,x,\Gamma)\right] \, ds|
$$
$$
\le \int_t^\infty \lambda e^{-\lambda s} | \int_E P(s,y,\Gamma) \, \nu(dy)-P(t,x,\Gamma) | \, ds
\le e^{-\lambda t}.
$$
Since $||\tilde{P}(t,x,\cdot)-q_\nu(\cdot)||_{TV}=2 \sup_{\Gamma\in {\cal E}} |\tilde{P}(t,x,\Gamma)-q_\nu(\Gamma)|$ \cite[Appendix]{Hernandez-Lerma:1996}, it follows that
$\{\tilde{X}(t),t\in [0,\infty)\}$ is exponentially ergodic with index $\lambda$ and $M(x)=2$.

\qed

\section{Moments and limits}

In this section, we let the initial distribution of the original process be the Dirac measure concentrated at $x\in E,$ where $E={\mathbb R}^n$. Then we write $\{\tilde{X}(t):t\ge 0\}$ as $\{(\tilde{X}_1(t),\dots,\tilde{X}_n(t)):t\ge 0\}.$ Similar notations are introduced for the process $\{X(t): t\ge 0\},$ too.
Now consider the $i$th component process $\{\tilde{X}_i(t),t\ge 0\}.$ From (\ref{eq:Ptilde}) we can obtain an expression for the moments
\begin{equation}\label{eq:moments}
E_x[\tilde X_i^k(t)] = e^{-\lambda t} E_x[X_i^k(t)]
+\int_E \int_0^t \lambda e^{-\lambda s} E_y[X_i^k(s)] \, ds \, \nu(dy),
\end{equation}
where and below in this section, we assume the involved interchange of the order of integrals is legal, which is the case,
for example, when
$$
\int_E \int_0^t \lambda e^{-\lambda s} E_y[|X_i^k(s)|] \, ds \, \nu(dy)<\infty.
$$
In turn, the equation (\ref{eq:moments}) helps to establish the following bound.

\begin{prop}
Let the $k$-th moment of the original process be exponentially bounded from the above in time, i.e.,
\begin{equation}\label{eq:exptimeboundo}
E_x[X^k_i(t)] \le c_{k,i}(x) e^{\eta_{k,i} t},
\end{equation}
where $c_{k,i}(\cdot)$ is a measurable $\nu$-integrable function, and $\eta_{k,i}<\lambda$ is a constant.
Then,
\begin{equation}\label{eq:exptimeboundo}
\limsup_{t \to \infty} E_x[\tilde X^k_i(t)] \le \frac{\bar{c}_{k,i}\lambda}{\lambda-\eta_{k,i}},
\end{equation}
where $\bar{c}_{k,i}=\int_E c_{k,i}(y) \nu(dy)$. If the $k$-th moment of the original
process is uniformly bounded from the above $(i.e., E_x[X^k_i(t)] \le c_{k,i}\in(-\infty,\infty))$, so is the $k$-th moment of the modified
process by the same bound $(i.e., E_x[\tilde X^k_i(t)] \le c_{k,i})$.
\end{prop}
{\bf Proof:} The equation (\ref{eq:moments}) yields
\begin{eqnarray*}
E_x[\tilde X^k_i(t)] &\le& e^{-\lambda t} c_{k,i}(x) e^{\eta_{k,i} t} +
\int_E \int_0^t \lambda e^{-\lambda s} c_{k,i}(y) e^{\eta_{k,i} s} \, ds \, \nu(dy)\\
&=&e^{-(\lambda-\eta_{k,i})t} c_{k,i}(x)+\frac{\lambda}{\lambda-\eta_{k,i}}(1-e^{-(\lambda-\eta_{k,i})t})\int_E c_{k,i}(y)\eta(dy),
\end{eqnarray*}
which implies the first statement of the proposition. The second statement of the proposition
also follows from the equation (\ref{eq:moments}).

\qed

Now let us investigate what happens when the parameter $\lambda$ goes to zero.
Since $\int_E P(s,y,\Gamma)\nu(dy)$ is a  bounded measurable function with respect to $s$,
according to \cite{F71}, we can conclude that the limiting discounting
is equivalent to the time averaging. This fact and Theorem 1 in \cite{Glynn:1994} lead to the following statement.

\begin{thm}\label{thm:Taub}
The existence of the limit $\lim_{\lambda \to 0} q_\nu(\Gamma)$ is equivalent to
the existence of the limit $\lim_{T \to \infty} 1/T\int_0^T \int_E P(s,y,\Gamma) \, \nu(dy) \, ds$.
If these limits exist, they are equal. In particular,
if the original process $\{X(t): t\in[0,\infty)\}$ is positive Harris recurrent, the limit $\lim_{\lambda \to 0} q_\nu(\Gamma)$
exists and is equal to the invariant probability measure of the original process.
\end{thm}

\begin{rem}
There could be cases when the limit $\lim_{\lambda \to 0} q_\nu(\cdot)$ exists for any
probability distribution $\nu(\cdot)$. However, similarly to the case of singularly perturbed
Markov processes \cite{AFH02}, the original process might not be ergodic.
\end{rem}

\section{Examples}

Let us illustrate the general results with two examples.

\subsection{Brownian motion with drift}

As the first example let us consider the Brownian motion with drift $\mu$
and variance coefficient $\sigma^2$ on the real line $E={\mathbb R}$ (see e.g., \cite{R96}).  The initial distribution of the original process is the Dirac measure concentrated at $x\in E.$
It can be described in the stochastic differential notation
$$
dX(t) = \mu dt + \sigma dW(t),
$$
where $W(t)$ is the standard Wiener process. The probability density function
of the Brownian process has a closed form
$$
p(t,0,z)=\frac{1}{\sqrt{2 \pi \sigma^2 t}} \exp\left( -\frac{(z-\mu t)^2}{2\sigma^2t}\right).
$$
Here we assume that the process starts from zero.
We observe that $p(t,0,z)$ does not converge to a proper probability density as $t$ goes to infinity.
If the modified process restarts also from zero, by formula (\ref{eq:Ptilde}) we have
$$
\tilde p(t,0,z)=\exp(-\lambda t)\frac{1}{\sqrt{2 \pi \sigma^2 t}} \exp\left( -\frac{(z-\mu t)^2}{2\sigma^2t}\right)
$$
\begin{equation}\label{eq:exampsol}
+\lambda \int_0^t \exp(-\lambda s)
\frac{1}{\sqrt{2 \pi \sigma^2 s}} \exp\left( -\frac{(z-\mu s)^2}{2\sigma^2s}\right) ds,
\end{equation}
which has the well defined limiting probability density function
$$
q_0(z)=\lambda \int_0^\infty \exp(-\lambda s)
\frac{1}{\sqrt{2 \pi \sigma^2 s}} \exp\left( -\frac{(z-\mu s)^2}{2\sigma^2s}\right) ds.
$$
Let us now consider that the original process starts with an arbitrarily fixed $x\in E$ and restarts according to the distribution $\nu(dy)$, which admits a finite second moment. Then we may calculate the first moment of the modified process by formula (\ref{eq:moments}) with
$E_x[X(t)]=x+\mu t$
$$
E_x[\tilde X(t)] = e^{-\lambda t} (x+\mu t)
+\int_E \int_0^t \lambda e^{-\lambda s} (y+\mu s) \, ds \, \nu(dy)
$$
$$
=e^{-\lambda t} (x+\mu t) + [1-e^{-\lambda t}] \int_E y \nu(dy)
+ [1-(1+\lambda t)e^{-\lambda t}]\frac{\mu}{\lambda}.
$$
Thus,
$$
E_x[\tilde X(t)] \to \int_E y \nu(dy) + \frac{\mu}{\lambda},
$$
as $t \to \infty$.

Similarly, a direct calculation gives
\begin{eqnarray*}
E_x[\tilde{X}^2(t)]\rightarrow \frac{\sigma^2}{\lambda}+\frac{2\mu^2}{\lambda^2}+ \int_E (\frac{2\mu y}{\lambda}+ y^2)\nu(dy), \mbox{ as }t\rightarrow \infty,
\end{eqnarray*}
and thus
\begin{eqnarray*}
Var_x[\tilde{X}^2(t)]\rightarrow \int_E y^2 \nu(dy)-\left(\int_E y\nu(dy)\right)^2+\frac{\sigma^2}{\lambda}+\frac{\mu^2}{\lambda^2}, \mbox{ as }t\rightarrow\infty.
\end{eqnarray*}

\subsection{Geometric Brownian motion}

As the second example we consider the geometric Brownian motion (see e.g., \cite{R96}), so that we take $E=[0,\infty).$
It can be described by the following stochastic differential equation (with the initial condition $P(X(0)\in dy)=1\{{x\in dy}\},$ where $x\in E$ is fixed)
$$
dX(t) = \mu X(t) dt + \sigma X(t) dW(t),
$$
where $W(t)$ is the standard Wiener process. The probability density function of the
geometric Brownian motion also has a closed form
\begin{equation}\label{eq:lognorm}
p(t,x,z)=
\frac{1}{\sqrt{2\pi}}\frac{1}{z\sigma\sqrt{t}}
\exp\left(-\frac{(\ln(z)-\ln(x)-(\mu-\sigma^2/2)t)^2}{2\sigma^2t}\right),
\end{equation}
which defines a log-normal distribution with mean
\begin{equation}\label{eq:lognormmean}
E_x[X(t)] = x e^{\mu t},
\end{equation}
and variance
\begin{equation}\label{eq:lognormvar}
Var_x[X(t)] = x^2 e^{2\mu t}(e^{\sigma^2 t}-1).
\end{equation}
Indeed, for any $k=1,2,\dots,$ it holds for this log-normal distribution that
\begin{eqnarray*}
E_x[X^k(t)]=x^k \exp(k
(\mu-\frac{\sigma^2}{2})t+\frac{k^2\sigma^2t}{2}).
\end{eqnarray*}

Next, for any fixed $k=1,2,\dots,$ we assume $\lambda-k(\mu-\frac{\sigma^2}{2})-\frac{1}{2}k^2\sigma^2\ne 0,$ and $\nu(dy)$ has a finite $k$th moment. By using formula (\ref{eq:moments}), we calculate the $k$th moment of the modified process
\begin{eqnarray*}
E_x[\tilde{X}^k(t)]&=& e^{-\lambda t} x^k e^{k(\mu-\frac{\sigma^2}{2})t+\frac{1}{2}k^2\sigma^2 t}+\int_E\int_0^t \lambda e^{-\lambda s} y^k e^{k(\mu-\frac{\sigma^2}{2})s+\frac{1}{2}k^2\sigma^2 s}ds \nu(dy)\\
&=&x^k e^{-t(\lambda- k(\mu-\frac{\sigma^2}{2})-\frac{1}{2}k^2\sigma^2)}+ \int_E \lambda y^k \int_0^t e^{-s (\lambda- k(\mu-\frac{\sigma^2}{2})-\frac{1}{2}k^2\sigma^2)}ds \nu(dy)\\
&=&x^k e^{-t(\lambda- k(\mu-\frac{\sigma^2}{2})-\frac{1}{2}k^2\sigma^2)}+  \lambda  \frac{1-e^{-t (\lambda- k(\mu-\frac{\sigma^2}{2})-\frac{1}{2}k^2\sigma^2)}}{\lambda- k(\mu-\frac{\sigma^2}{2})-\frac{1}{2}k^2\sigma^2}  \int_E y^k  \nu(dy).
\end{eqnarray*}
Thus, when $\lambda > k(\mu-\frac{\sigma^2}{2}) + \frac{1}{2}k^2\sigma^2> 0,$ we have
\begin{eqnarray*}
E_x[\tilde X^k(t)] \to \frac{\lambda}{{\lambda- k(\mu-\frac{\sigma^2}{2})-\frac{1}{2}k^2\sigma^2}} \int_E y^k\nu(dy)
\end{eqnarray*}
as $t \to \infty$.

\newpage
\tableofcontents

\end{document}